\documentclass{article}

\usepackage{amsmath,amssymb,amsthm}
\usepackage{graphicx,slashbox}
\usepackage{subfigure}
\usepackage{textcomp}
\usepackage{hyperref}
\usepackage{a4wide}

\newtheorem{Theorem}{Theorem}
\theoremstyle{Definition}
\newtheorem{Definition}[Theorem]{Definition}

\def\a{\alpha}
\def\t{\tau}
\def\l{\lambda}

\def\RIT{{_tI_T^{1-\a}}}
\def\LCD{{^C_aD_t^\a}}
\def\RDT{{_tD_T^\a}}

\def\LDa{{_aD_t^\a}}

\newenvironment{keywords}{\begin{center}
\begin{minipage}[c]{13.4cm} {\bf Keywords:}} {\end{minipage}
\end{center}}

\newenvironment{msc}{\begin{center}
\begin{minipage}[c]{13.4cm} {\bf MSC 2010:}} {\end{minipage}
\end{center}}

\begin{document}

\title{A numerical scheme to solve fractional\\
optimal control problems\footnote{Part of first author's Ph.D.,
which is carried out at the University of Aveiro under the
\emph{Doctoral Program in Mathematics and Applications}
(PDMA) of Universities of Aveiro and Minho. Submitted 13-Apr-2013;
accepted after minor revision 22-May-2013; 
Conference Papers in Mathematics, Volume 2013, Article ID 165298, 10 pages.
\url{http://dx.doi.org/10.1155/2013/165298}}}

\author{Shakoor Pooseh\\
\texttt{spooseh@ua.pt}
\and Ricardo Almeida\\
\texttt{ricardo.almeida@ua.pt}
\and Delfim F. M. Torres\\
\texttt{delfim@ua.pt}}

\date{CIDMA -- Center for Research and Development in Mathematics and Applications,\\
Department of Mathematics, University of Aveiro, 3810-193 Aveiro, Portugal}

\maketitle


\begin{abstract}
We review recent results obtained to solve fractional order optimal control problems
with free terminal time and a dynamic constraint involving integer and fractional order derivatives.
Some particular cases are studied in detail. A numerical scheme is given, based on expansion formulas
for the fractional derivative. The efficiency of the method is illustrated through examples.
\end{abstract}

\begin{msc}
26A33, 33F05, 49K15.
\end{msc}



\section{Introduction}

In a letter dated September 30, 1695 l'H\^{o}pital posed the question to Leibniz:
what would be the derivative of order $\alpha=1/2$? Leibniz's response was:
``an apparent paradox, from which one day useful consequences will be drawn.''
In these words fractional calculus was born. In 1730, based on the formula
$$
\frac{d^n x^m}{dx^n}= m(m-1)\cdots (m-n+1)x^{m-n}
=\frac{\Gamma(m+1)}{\Gamma(m-n+1)}x^{m-n},
$$
Euler suggested to use this relationship also for real values of $n$.
Taking $m= 1$ and $n =1/2$, he obtained
$$
\frac{d^{1/2} x}{dx^{1/2}}=\sqrt{\frac{4x}{\pi}}.
$$
Since then, many different approaches have been carried out,
trying to find proper definitions for what should
be a derivative and an integral of real order.
Starting with  Cauchy's formula for an $n$-fold integral,
$$
\int_a^t d\tau_1\int_a^{\tau_1}d\tau_2\cdots \int_a^{\tau_{n-1}} x(\tau_n)d\tau_n
= \frac{1}{(n-1)!}\int_a^t (t-\tau)^{n-1}x(\tau) d\tau,
$$
Riemann  defined fractional integration as
$$
{_aI_t^{\alpha}}x(t) = \frac{1}{\Gamma(\alpha)}\int_a^t (t-\tau)^{\alpha-1}x(\tau) d\tau .
$$
This is nowadays the most common definition for fractional integral. We remark that when
the order $\alpha$ is an integer, then the fractional integral becomes a multiple integral,
recovering by this way the classical case.

We begin with some basic definitions and properties
about fractional operators \cite{kilbas,Podlubny}.
To avoid too many details, we omit here the conditions
that ensure the existence of such fractional operators
and the assumptions in which the results given below hold.
For an introduction to the fractional variational calculus
we refer the reader to \cite{b:agnieszka:delfim}.

\begin{Definition}
\label{def:oper}
Let $x:[a,b]\to\mathbb{R}$ be a function, $\alpha>0$ a real,
and $n=[\alpha]+1$, where $[\cdot]$ denotes the integer part function.
The left and right Riemann--Liouville fractional integrals are defined, 
respectively, by
\begin{gather}
{_aI_t^{\alpha}}x(t)
= \frac{1}{\Gamma(\alpha)}\int_a^t (t-\tau)^{\alpha-1}x(\tau) d\tau,
\tag{left RLFI}\\
{_tI_b^{\alpha}}x(t)
= \frac{1}{\Gamma(\alpha)}\int_t^b (\tau-t)^{\alpha-1}x(\tau) d\tau.
\tag{right RLFI}
\end{gather}
The left and right Riemann--Liouville fractional derivatives 
are defined, respectively, by
\begin{gather}
{_aD_t^{\alpha}}x(t)
= \frac{1}{\Gamma(n-\alpha)}\frac{d^n}{dt^n}
\int_a^t (t-\tau)^{n-\alpha-1}x(\tau) d\tau,
\tag{left RLFD}\\
{_tD_b^{\alpha}}x(t)
= \frac{(-1)^n}{\Gamma(n-\alpha)}\frac{d^n}{dx^n}
\int_t^b (\tau-t)^{n-\alpha-1}x(\tau) d\tau.
\tag{right RLFD}
\end{gather}
The left and right Caputo fractional derivatives 
are defined, respectively, by
\begin{gather}
{_a^CD_t^{\alpha}}x(t)
= \frac{1}{\Gamma(n-\alpha)}
\int_a^t (t-\tau)^{n-\alpha-1}x^{(n)}(\tau) d\tau,
\tag{left CFD}\\
{_t^CD_b^{\alpha}}x(t)
= \frac{(-1)^n}{\Gamma(n-\alpha)}
\int_t^b (\tau-t)^{n-\alpha-1}x^{(n)}(\tau) d\tau.
\tag{right CFD}
\end{gather}
\end{Definition}

We remark that if $\alpha=n$ in Definition~\ref{def:oper}, then we have the usual operators:
\begin{equation}
\begin{split}
{_aI_t^{n}}x(t)
&= \int_a^t d\tau_1 \int_a^{\tau_1} d\tau_2 \cdots \int_a^{\tau_{n-1}}x(\tau_n) d\tau_n,\\
{_tI_b^{n}}x(t)
&= \int_t^b d\tau_1 \int_{\tau_1}^b d\tau_2 \cdots \int_{\tau_{n-1}}^b x(\tau_n) d\tau_n,\\
{_aD_t^{n}}x(t)
&= x^{(n)}(t),\\
{_tD_b^{n}}x(t)
&= (-1)^n x^{(n)}(t),\\
{_a^CD_t^{n}}x(t)
&= x^{(n)}(t),\\
{_t^CD_b^{\alpha}}x(t)
&= (-1)^n x^{(n)}(t) .
\end{split}
\end{equation}

Some basic properties are useful, namely, a relationship between the Riemann--Liouville
and the Caputo fractional derivatives and a fractional integration by parts formula.

\begin{Theorem}
The following conditions hold:
\begin{enumerate}
\item $\displaystyle{_a^CD_t^\alpha}x(t)
={_aD_t^\alpha}x(t)-\sum_{k=0}^{n-1}\frac{x^{(k)}(a)}{\Gamma(k-\alpha+1)}(t-a)^{k-\alpha}$,

\item $\displaystyle {_aI_t^\alpha}{_aI_t^\beta}x(t)={_aI_t^{\alpha+\beta}}x(t)$,

\item $\displaystyle {_a^CD_t^\alpha}{_aI_t^\alpha}x(t)=x(t)$,

\item $\displaystyle {_aI_t^\alpha}{_a^CD_t^\alpha}x(t)=x(t)
-\sum_{k=0}^{n-1}\frac{x^{(k)}(a)}{k!}(t-a)^{k}$,

\item $\displaystyle \int_{a}^{b}y(t)\cdot {_a^C D_t^\alpha}x(t)dt
=\int_a^b x(t)\cdot {_t D_b^\alpha} y(t)dt
+\sum_{j=0}^{n-1}\left[{_tD_b^{\alpha+j-n}}y(t) \cdot {_tD_b^{n-1-j}} x(t)\right]_a^b$.
\end{enumerate}
\end{Theorem}

For numerical purposes, one of the most common procedure is to replace
the fractional operators by a series that involves integer derivatives only.
The usual one is given by
$$
{_aD_t^{\alpha}} x(t)=\sum_{n=0}^\infty \binom{\a}{n}
\frac{(t-a)^{n-\a}}{\Gamma(n+1-\a)}x^{(n)}(t),
$$
where
$$
\binom{\a}{n}=\frac{(-1)^{n-1}\a\Gamma(n-\a)}{\Gamma(1-\a)\Gamma(n+1)}.
$$
Although very simple to use, it is easy to conclude that in order to have
a small error when we approximate ${_aD_t^{\alpha}} x$ by a finite sum up to order $N$,
we need to consider a large value for $N$, i.e., we need to consider
the set of admissible functions to be $C^N[a,b]$ which is an important
restriction of the set of the space of functions.
Recently, in \cite{Atanackovic}, a new expansion formula is given,
with the advantage that we only need the first derivative:
\begin{equation}
\label{2}
\LDa x(t)=A(\a)(t-a)^{-\a}x(t)+B(\a)(t-a)^{1-\a}\dot{x}(t)
-\sum_{p=2}^{\infty}C(\a,p)(t-a)^{1-p-\a}V_p(t),
\end{equation}
where
$V_p(t)$ is the solution of the system
\begin{equation}
\label{eq:new+}
\begin{cases}
\dot{V}_p(t)=(1-p)(t-a)^{p-2}x(t),\\
V_p(a)=0,
\end{cases}
\end{equation}
for $p=2,3,\ldots,$ and $A$, $B$ and $C$ are given by
\begin{equation}
\label{eq:new-}
\begin{split}
A(\a) &= \frac{1}{\Gamma(1-\a)}
\left[1+\sum_{p=2}^{\infty}\frac{\Gamma(p-1+\a)}{\Gamma(\a)(p-1)!}\right],\\
B(\a) &= \frac{1}{\Gamma(2-\a)}
\left[1+\sum_{p=1}^{\infty}\frac{\Gamma(p-1+\a)}{\Gamma(\a-1)p!}\right],\\
C(\a,p) &= \frac{1}{\Gamma(2-\a)\Gamma(\a-1)}\frac{\Gamma(p-1+\a)}{(p-1)!}.
\end{split}
\end{equation}

We mention the recent papers \cite{Almeida1,Almeida3,Almeida2}, where similar results
are proven for fractional integrals and for other types of fractional operators.


\section{Necessary and sufficient optimality conditions}

Let $\a\in(0,1)$, and let $L,f: [a,+\infty[\times \mathbb{R}^2\to\mathbb{R}$
be two differentiable functions and $\phi:[a,+\infty[\times \mathbb{R}\to\mathbb{R}$
a differentiable function. The fundamental problem, as studied in \cite{Almeida4}, is the following:
\begin{equation}
\label{Opt:func}
\mathrm{minimize} \quad
J(x,u,T)=\int_a^T L(t,x(t),u(t))\,dt+\phi(T,x(T))
\end{equation}
subject to the dynamic control system
\begin{equation}
\label{Opt:dynamic}
M \dot{x}(t) + N\LCD x(t) = f\left(t,x(t),u(t)\right)
\end{equation}
and the initial condition
\begin{equation}
\label{Opt:bound}
x(a)=x_a,
\end{equation}
with $(M,N)\not=(0,0)$ and $x_a$ a fixed real number.
Thus, we are not only interested in finding the optimal
state function $x$ and the optimal control $u$,
but also the optimal time $T$.

\begin{Theorem}
\label{teo}
If $(x,u,T)$ is a minimizer of \eqref{Opt:func} under the dynamic constraint
\eqref{Opt:dynamic} and the boundary condition \eqref{Opt:bound},
then there exists a function $\l$ for which the triplet $(x,u,\l)$ satisfies
\begin{itemize}
\item[(i)] the \emph{Hamiltonian system}
\begin{equation*}
\begin{cases}
M \dot{\l}(t) - N\RDT\l(t) = -\frac{\partial H}{\partial x}\left(t,x(t),u(t),\l(t)\right)\\
M \dot{x}(t)+N\LCD x(t) = \frac{\partial H}{\partial \l}\left(t,x(t),u(t),\l(t)\right)
\end{cases}
\end{equation*}
for all $t\in[a,T]$;
\item[(ii)] the \emph{stationary condition}
\begin{equation*}
\frac{\partial H}{\partial u}(t,x(t),u(t),\l(t))=0,
\quad \forall \, t\in[a,T];
\end{equation*}
\item[(iii)] the \emph{transversality conditions}
\begin{equation}
\label{t1}
\left[H(t,x(t),u(t),\l(t))-N\l(t)\LCD x(t)+N\dot{x}(t)\RIT\l(t)
+\frac{\partial \phi}{\partial t}(t,x(t))\right]_{t=T}=0,
\end{equation}
\begin{equation}
\label{t2}
\left[M \l(t) +N\RIT\l(t)-\frac{\partial \phi}{\partial x}(t,x(t))\right]_{t=T}=0;
\end{equation}
\end{itemize}
where the Hamiltonian $H$ is defined by
$$
H(t,x,u,\l)=L(t,x,u)+\l f(t,x,u).
$$
\end{Theorem}

This theorem states the general condition that the optimal solution $(x,u,T)$ must fulfill.
Next, depending on extra conditions imposed over the final time $T$ or in $x(T)$,
new transversality conditions are obtained.

\begin{Theorem}
\label{teo1}
Let $(x,u)$ be a minimizer of \eqref{Opt:func} under the dynamic constraint
\eqref{Opt:dynamic} and the boundary condition \eqref{Opt:bound}.
\begin{itemize}

\item[(i)] If $T$ is fixed and $x(T)$ is free, then Theorem~\ref{teo}
holds with the transversality conditions
\eqref{t1} and \eqref{t2} replaced by
$$
\left[M \l(t)+N\RIT\l(t)
-\frac{\partial \phi}{\partial x}(t,x(t))\right]_{t=T}=0.
$$

\item[(ii)] If $x(T)$ is fixed and $T$ is free,
then Theorem~\ref{teo} holds with the transversality
conditions \eqref{t1} and \eqref{t2}  replaced by
$$
\left[H(t,x(t),u(t),\l(t))-N\l(t)\LCD x(t)+N\dot{x}(t)\RIT\l(t)
+\frac{\partial \phi}{\partial t}(t,x(t))\right]_{t=T}=0.
$$

\item[(iii)] If $T$ and $x(T)$ are fixed,
then Theorem~\ref{teo} holds with no transversality conditions.

\item[(iv)] If the terminal point $x(T)$ belongs to a fixed curve, i.e.,
$x(T)=\gamma(T)$ for some differentiable curve $\gamma$, then
Theorem~\ref{teo} holds with the transversality conditions
\eqref{t1} and \eqref{t2} replaced by
\begin{multline*}
\Biggl[H(t,x(t),u(t),\l(t))-N\l(t)\LCD x(t)
+N\dot{x}(t)\RIT\l(t)+\frac{\partial \phi}{\partial t}(t,x(t))\\
-\dot{\gamma}(t)\left(M \l(t)+N\RIT\l(t)
-\frac{\partial \phi}{\partial x}(t,x(t))\right)\Biggr]_{t=T}=0.
\end{multline*}

\item[(v)] If $T$ is fixed and $x(T)\geq K$ for some fixed $K\in\mathbb{R}$,
then Theorem~\ref{teo} holds with the transversality conditions
\eqref{t1} and \eqref{t2} replaced by
\begin{gather*}
\left[M \l(t)+N\RIT\l(t)-\frac{\partial \phi}{\partial x}(t,x(t))\right]_{t=T}\leq 0,\\
(x(T)-K)\left[M \l(t)+N\RIT\l(t)-\frac{\partial \phi}{\partial x}(t,x(t))\right]_{t=T}=0.
\end{gather*}
\end{itemize}
\end{Theorem}

Numerically, by using approximation \eqref{2} up to order $K$,
we can transform the original problem into the following classical optimal control problem:
\begin{equation*}
\mathrm{minimize} \quad
\tilde{J}(x,u,T)=\int_a^TL(t,x(t),u(t))\,dt+\phi(T,x(T))
\end{equation*}
subject to the dynamic constraints
\begin{equation*}
\begin{cases}
\dot{x}(t)=\displaystyle \frac{f(t,x(t),u(t))
-NA(t-a)^{-\a}x(t)+\sum_{p=2}^KNC_p(t-a)^{1-p-\a}V_p(t)}{M+NB(t-a)^{1-\a}}\\[0.25cm]
\dot{V}_p(t)=(1-p)(t-a)^{p-2}x(t), \quad p=2,\ldots,K
\end{cases}
\end{equation*}
and the initial conditions
\begin{equation}
\label{App:bound}
\begin{cases}
x(a)=x_a,\\
V_p(a)=0, \quad p=2,\ldots,K.
\end{cases}
\end{equation}

Theorem~\ref{teo} can be generalized in the following way.
Observe that we have two initial points for the problem,
one for the fractional derivative and a second one for
the integral of the functional. We now consider a more general approach,
where the initial time for the integral is greater than the initial time
of the fractional derivative. We impose a boundary condition on $t=A$,
but similar conditions could be obtained if we considered conditions
at $t=a$ instead. The problem is formulated as follows. Let $\a\in(0,1)$, and let
$L,f: [a,+\infty[\times \mathbb{R}^2\to\mathbb{R}$ be two differentiable functions,
$\phi:[a,+\infty[\times \mathbb{R}\to\mathbb{R}$ a differentiable function, and $A>a$ a real.
We wish to
\begin{equation}
\label{Opt:func2}
\mathrm{minimize} \quad
J(x,u,T)=\int_A^T L(t,x(t),u(t))\,dt+\phi(T,x(T))
\end{equation}
subject to
\begin{equation}
\label{Opt:dynamic2}
M \dot{x}(t) + N\LCD x(t) = f\left(t,x(t),u(t)\right)
\end{equation}
and
\begin{equation}
\label{Opt:bound2}
x(A)=x_A,
\end{equation}
with $(M,N)\not=(0,0)$ and $x_A$ being a fixed real number.

\begin{Theorem}
\label{teo2}
If $(x,u,T)$ is a minimizer of \eqref{Opt:func2} under the dynamic constraint
\eqref{Opt:dynamic2} and the boundary condition \eqref{Opt:bound2},
then there exists a function $\l$ for which the triplet $(x,u,\l)$ satisfies
\begin{itemize}
\item[(i)] the \emph{Hamiltonian system}
\begin{equation*}
\begin{cases}
M \dot{\l}(t)-N\RDT\l(t)
= - \displaystyle \frac{\partial H}{\partial x}(t,x(t),u(t),\l(t))\\[0.25cm]
M \dot{x}(t)+N\LCD x(t)
= \displaystyle \frac{\partial H}{\partial \l}(t,x(t),u(t),\l(t))
\end{cases}
\end{equation*}
for all $t\in[A,T]$ and
$$
\RDT\l(t)-{_tD^\a_A}\l(t)=0
$$
for all $t\in[a,A]$;

\item[(ii)] the \emph{stationary condition}
$$
\frac{\partial H}{\partial u}(t,x(t),u(t),\l(t))=0,
\quad \forall \, t\in[A,T];
$$
\item[(iii)] the \emph{transversality conditions}
\begin{gather*}
\left[H(t,x(t),u(t),\l(t))-N\l(t)\LCD x(t)+N\dot{x}(t)\RIT\l(t)
+\frac{\partial \phi}{\partial t}(t,x(t))\right]_{t=T}=0,\\
\left[M \l(t)+N\RIT\l(t)-\frac{\partial \phi}{\partial x}(t,x(t))\right]_{t=T}=0,\\
\left[{_tI^{1-\a}_T}\l(t)-{_tI^{1-\a}_A}\l(t)\right]_{t=a}=0;
\end{gather*}
\end{itemize}
where the Hamiltonian $H$ is defined by
$$
H(t,x,u,\l)=L(t,x,u)+\l f(t,x,u).
$$
\end{Theorem}

We remark that when $A=a$, Theorem~\ref{teo2} reduces to Theorem~\ref{teo}.

Under some additional conditions, namely, convexity conditions over $L$, $f$, and $\phi$,
Theorem~\ref{teo} provides also sufficient conditions to ensure optimal solutions.
The result is given in the next theorem.

\begin{Theorem}
Let $\left(\overline{x}, \overline{u}, \overline{\l}\right)$
be a triplet satisfying the necessary conditions of Theorem~\ref{teo}.
Moreover, assume that
\begin{enumerate}
\item $L$ and $f$ are convex on $x$ and $u$, and $\phi$ is convex in $x$;
\item $T$ is fixed;
\item $\overline{\l}(t)\geq 0$ for all $t \in [a,T]$ or $f$ is linear in $x$ and $u$.
\end{enumerate}
Then $\left(\overline x,\overline u\right)$ is an optimal solution to problem
\eqref{Opt:func}--\eqref{Opt:bound}.
\end{Theorem}


\section{Numerical treatment}

So far, we have provided a theoretical approach to fractional optimal control problems,
which involves solving fractional differential equations. As it is known, solving such
equations is in most cases impossible to do, and numerical methods are used to find
approximated solutions for the problem (see, e.g., \cite{Almeida0,Almeida5}).
We describe next, briefly, how formula \eqref{2} is deduced
and generalized for arbitrary size expansions.

Let $x\in C^2[a,b]$. Using integration by parts two times, we deduce that
$$
{_aD^\alpha_t} x(t)=\frac{x(a)}{\Gamma(1-\a)}(t-a)^{-\a}
+\frac{\dot{x}(a)}{\Gamma(2-\a)}(t-a)^{1-\a}+\frac{(t-a)^{1-\a}}{\Gamma(2-\a)}
\int_a^t \left(1-\frac{\t-a}{t-a}\right)^{1-\a}\ddot{x}(\t)d\t.
$$
By the binomial formula, we can rewrite the fractional derivative as
\begin{multline*}
{_aD^\alpha_t} x(t)=\frac{x(a)}{\Gamma(1-\a)}(t-a)^{-\a}
+\frac{(t-a)^{1-\a}}{\Gamma(2-\a)}\dot{x}(a)\\
+\frac{(t-a)^{1-\a}}{\Gamma(2-\a)}\int_a^t
\left(\sum_{p=0}^{\infty}\frac{\Gamma(p-1+\a)}{\Gamma(\a-1)p!}\left(
\frac{\t-a}{t-a}\right)^p\right)\ddot{x}(\t)d\t.
\end{multline*}
Further integration by parts gives
$$
{_aD^\alpha_t} x(t)=A(\a)(t-a)^{-\a}x(t)+B(\a)(t-a)^{1-\a}\dot{x}(t)
-\sum_{p=2}^{\infty}C(\a,p)(t-a)^{1-p-\a}V_p(t),
$$
where $V_p(t)$, $A(\a)$, $B(\a)$ and $C(\a,p)$ 
are given in \eqref{eq:new+}--\eqref{eq:new-}.
Following similar calculations, we can deduce the next theorem.

\begin{Theorem}
Fix $n \in \mathbb{N}$ and let $x\in C^n[a,b]$. Then,
\begin{multline*}
{_aD^\alpha_t}x(t)=\frac{1}{\Gamma(1-\a)}(t-a)^{-\a}x(t)
+\sum_{i=1}^{n-1}A(\a,i)(t-a)^{i-\a}x^{(i)}(t)\\
+\sum_{p=n}^{\infty}\left[\frac{-\Gamma(p-n+1+\a)}{
\Gamma(-\a)\Gamma(1+\a)\left(p-n+1\right)!}(t-a)^{-\a}x(t)
+ B(\a,p)(t-a)^{n-1-p-\a}V_p(t)\right],
\end{multline*}
where
\begin{equation*}
\begin{split}
A(\a,i)&=\frac{1}{\Gamma(i+1-\a)}\left[1+\sum_{p=n-i}^{\infty}
\frac{\Gamma(p-n+1+\a)}{\Gamma(\a-i)(p-n+i+1)!}\right],\\
B(\a,p)&=\frac{\Gamma(p-n+1+\a)}{\Gamma(-\a)\Gamma(1+\a)(p-n+1)!},\\
V_p(t) &=(p-n+1)\int_a^t (\t-a)^{p-n}x(\t)d\t.
\end{split}
\end{equation*}
\end{Theorem}

The idea is to replace the fractional derivative with such expansions
and to consider finite sums only. When we use the approximation
$$
{_aD^\alpha_t}x(t) \approx \sum_{i=0}^{n-1} A(\a,i,N)(t-a)^{i-\a}x^{(i)}(t)
+\sum_{p=n}^{N} B(\a,p)(t-a)^{n-1-p-\a}V_p(t),
$$
the error is bounded by
$$
|E_{tr}(t)|\leq L_n\frac{\mathrm{e}^{(n-1-\a)^2
+n-1-\a}}{\Gamma(n-\a)(n-1-\a)N^{n-1-\a}}(t-a)^{n-\a},
$$
where
$$
L_n=\displaystyle\max_{\t \in [a,t]}\left|x^{(n)}(\t)\right|.
$$

\begin{figure}[ht]
\begin{center}
\subfigure[$x_1(t)=t^4$]{\label{plot1}\includegraphics[scale=0.50]{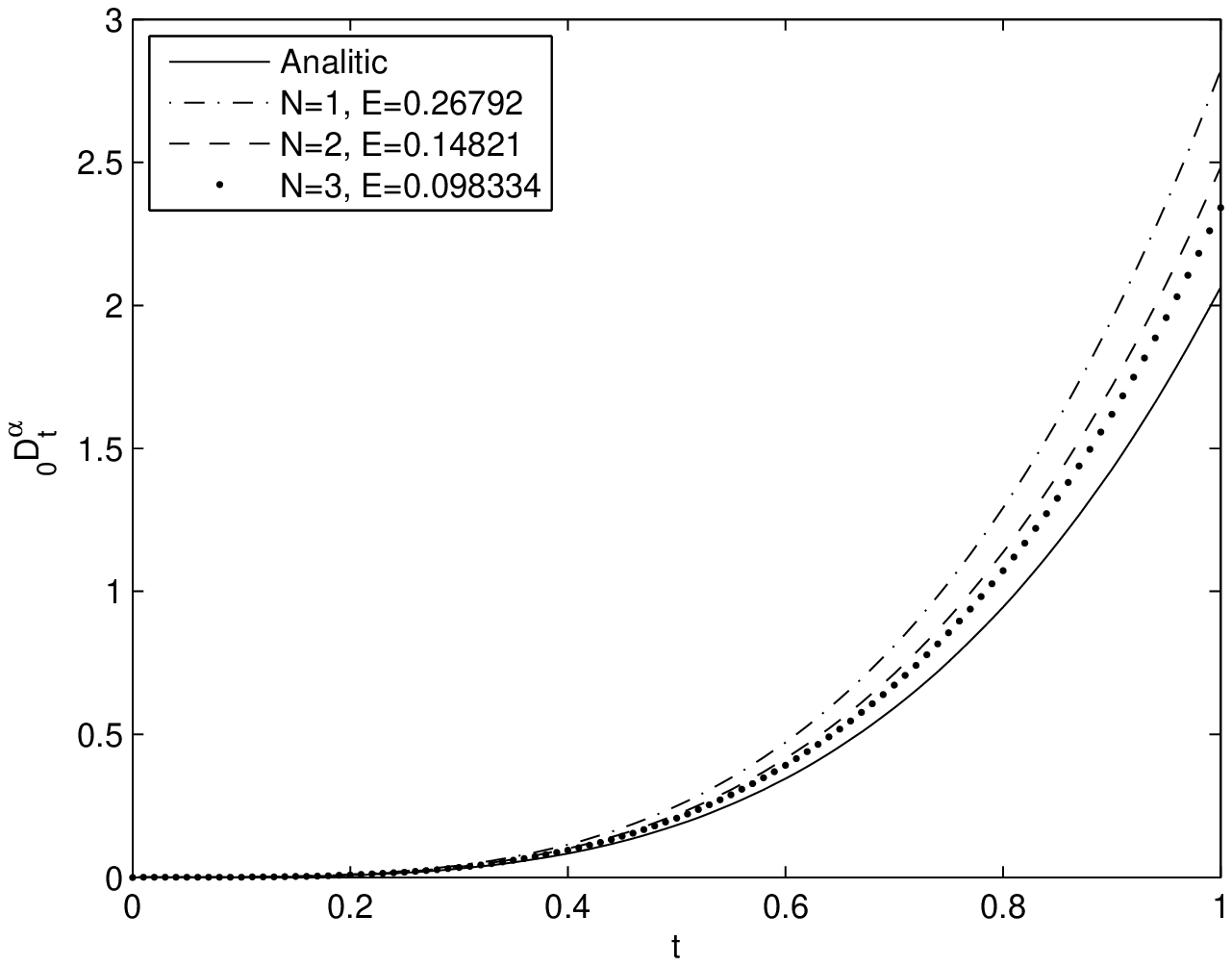}}
\subfigure[$x_2(t) = e^{2t}$]{\label{plot2}\includegraphics[scale=0.50]{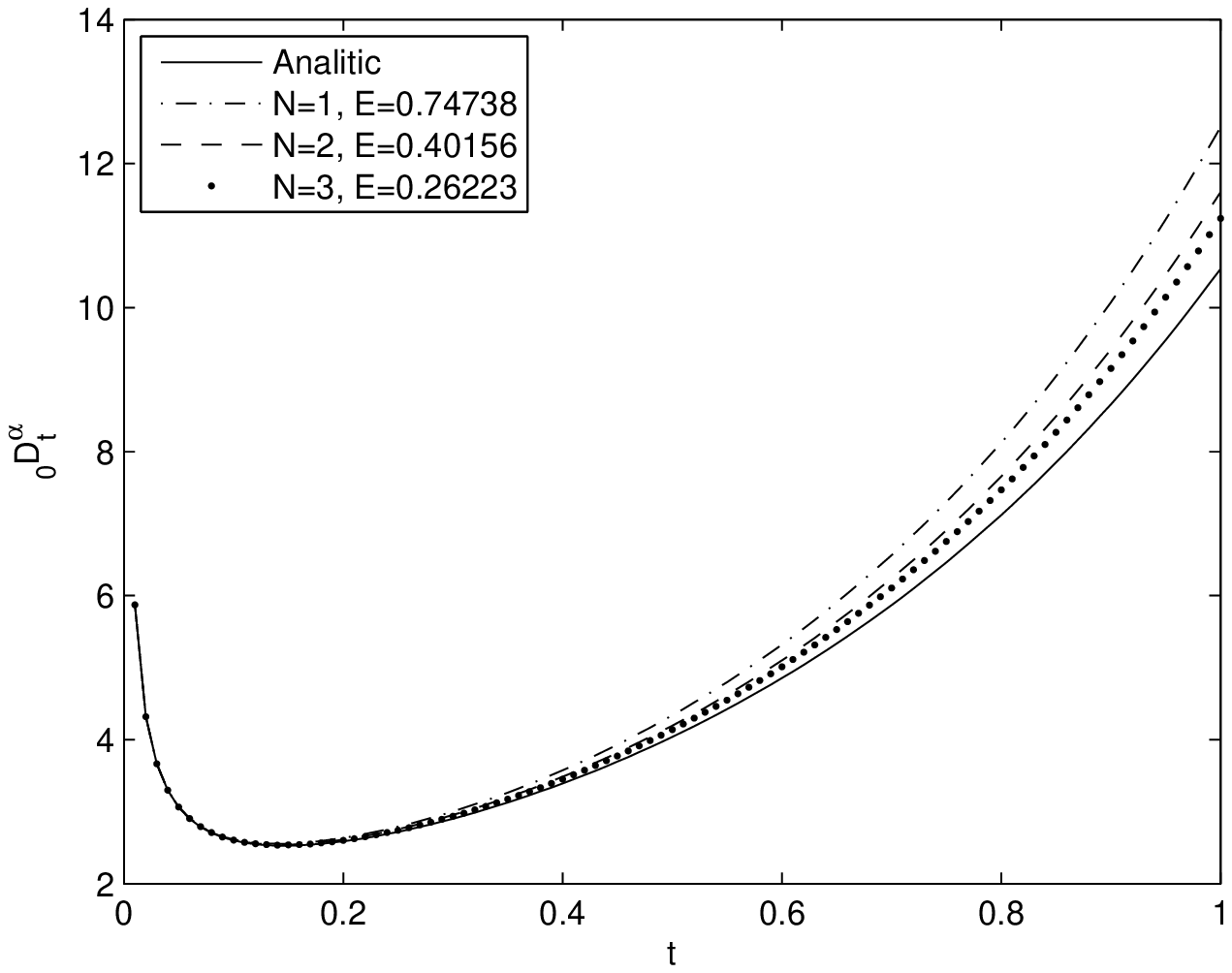}}
\end{center}
\caption{Analytic versus numerical approximation for a fixed $n$ ($n = 2$).}
\end{figure}

To see the accuracy of the method, we exemplify it by considering
some functions and compare the exact expression of the fractional
derivative with the approximated one. To start, consider $x_1(t)=t^4$
and  $x_2(t)=e^{2t}$ and expansions with $n=2$ and different values for $N$.
The result is exemplified in Figures~\ref{plot1} and \ref{plot2}.

\begin{figure}[ht]
\begin{center}
\subfigure[$x_1(t)=t^4$]{\label{plot3}\includegraphics[scale=0.50]{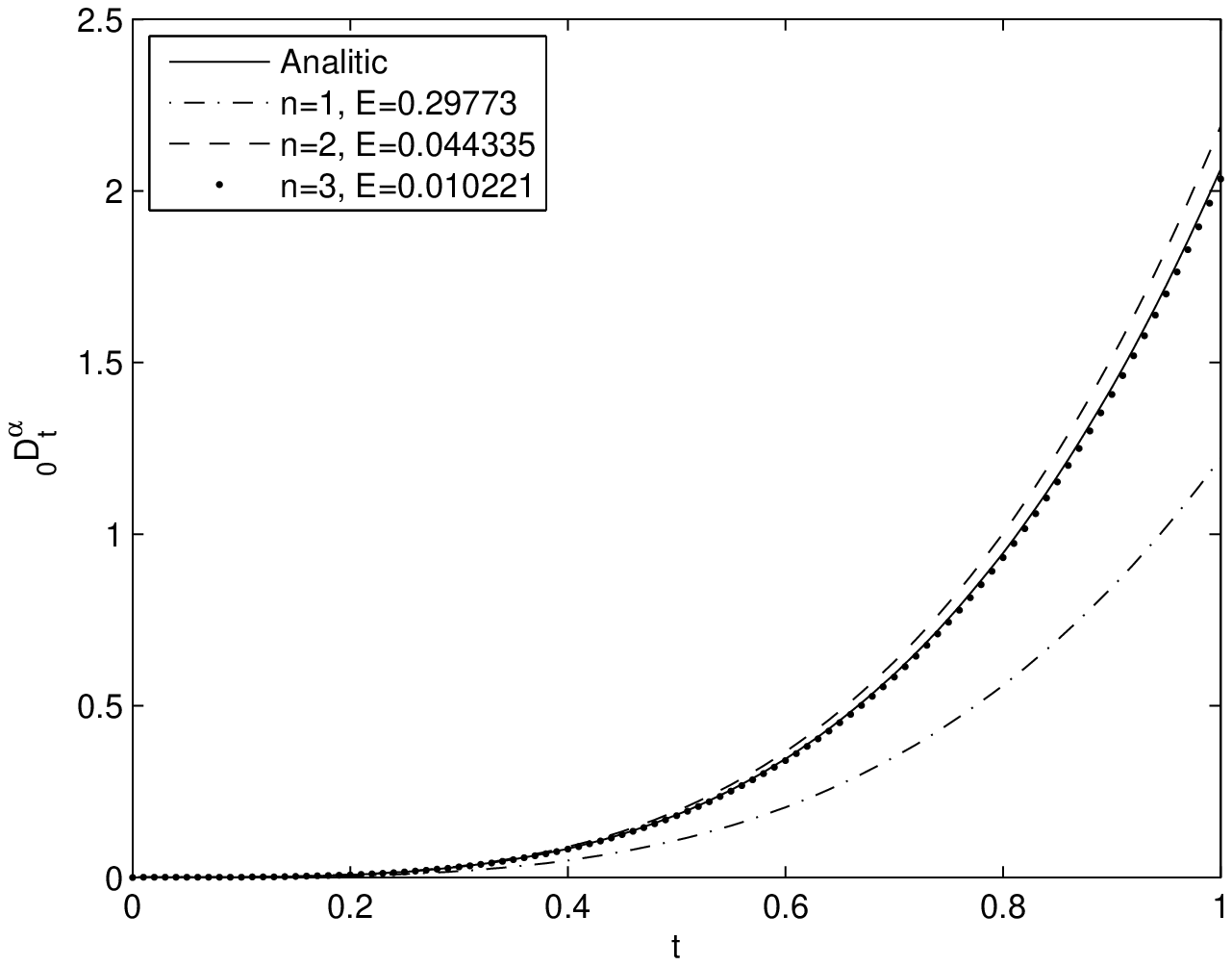}}
\subfigure[$x_2(t) = e^{2t}$]{\label{plot4}\includegraphics[scale=0.50]{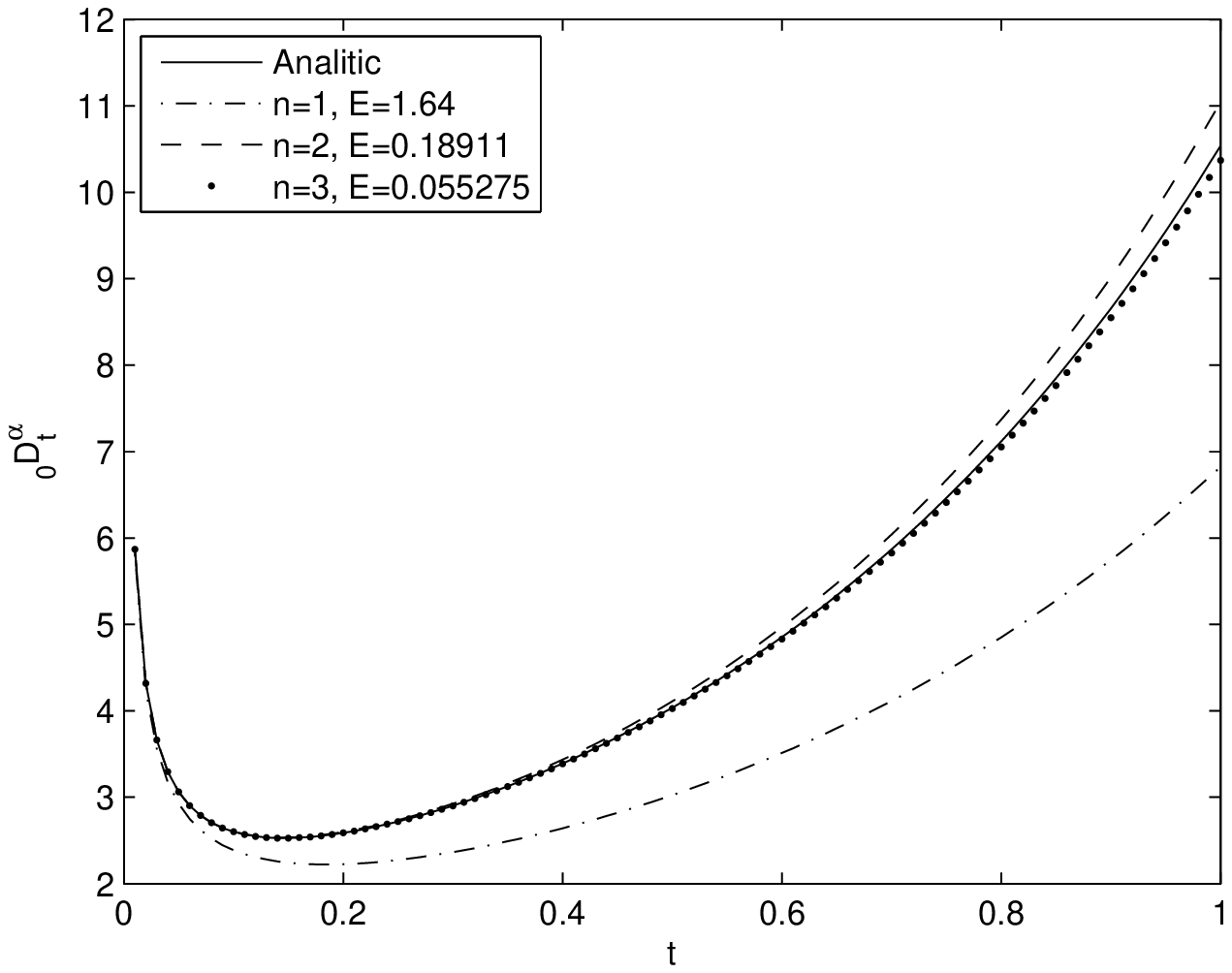}}
\end{center}
\caption{Analytic versus numerical approximation for a fixed $N$ ($N=6$).}
\end{figure}

A different approach is to consider a fixed $N$ and vary the size of the expansion,
i.e., to consider different values for $n$. For the same functions $x_1$ and $x_2$,
with $N=6$, the results are shown in Figures~\ref{plot3} and \ref{plot4}.


\section{Examples}

\begin{figure}[ht]
\begin{center}
\subfigure[$x(t), N=3$]{\includegraphics[scale=0.37]{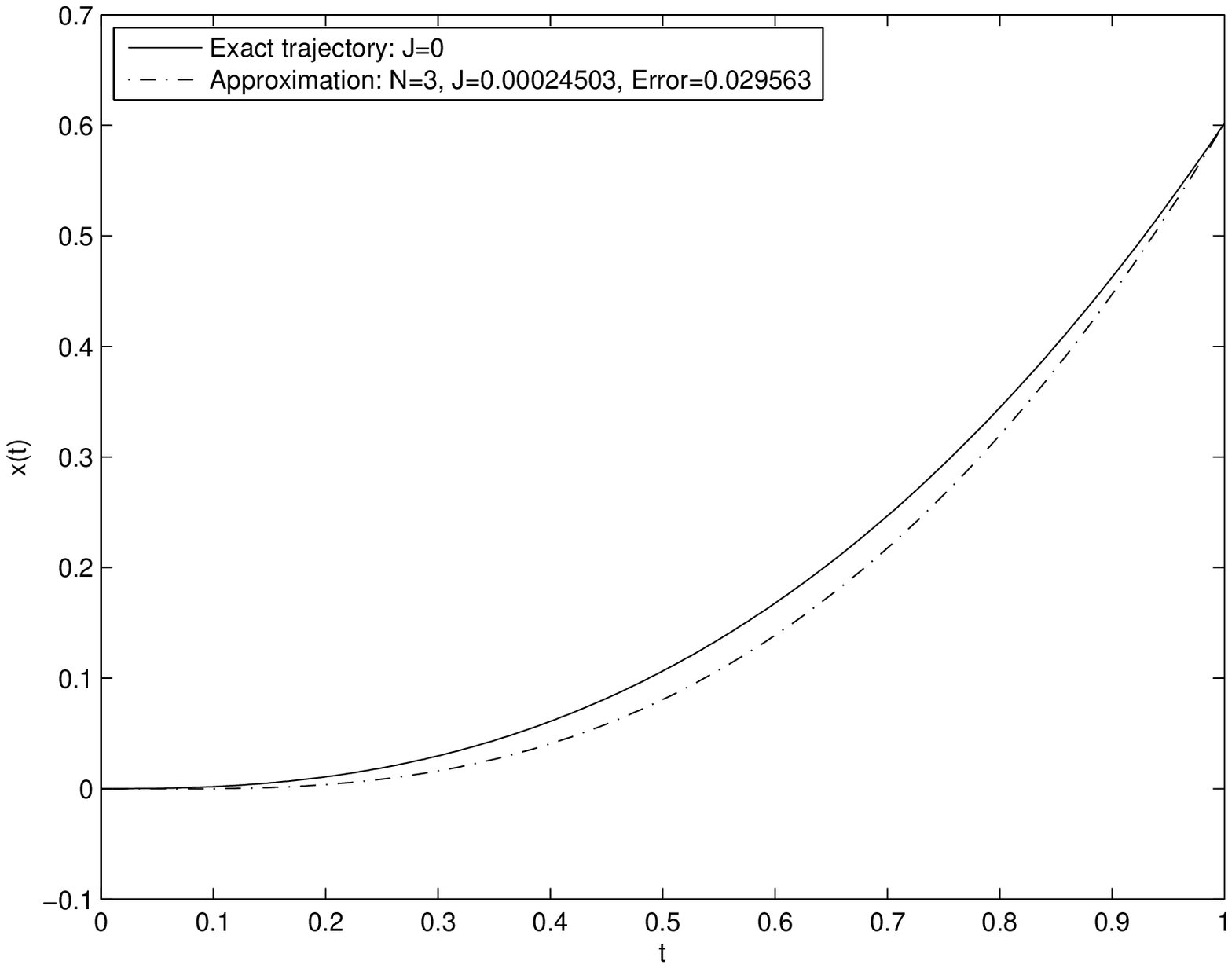}}
\subfigure[$u(t), N=3$]{\includegraphics[scale=0.37]{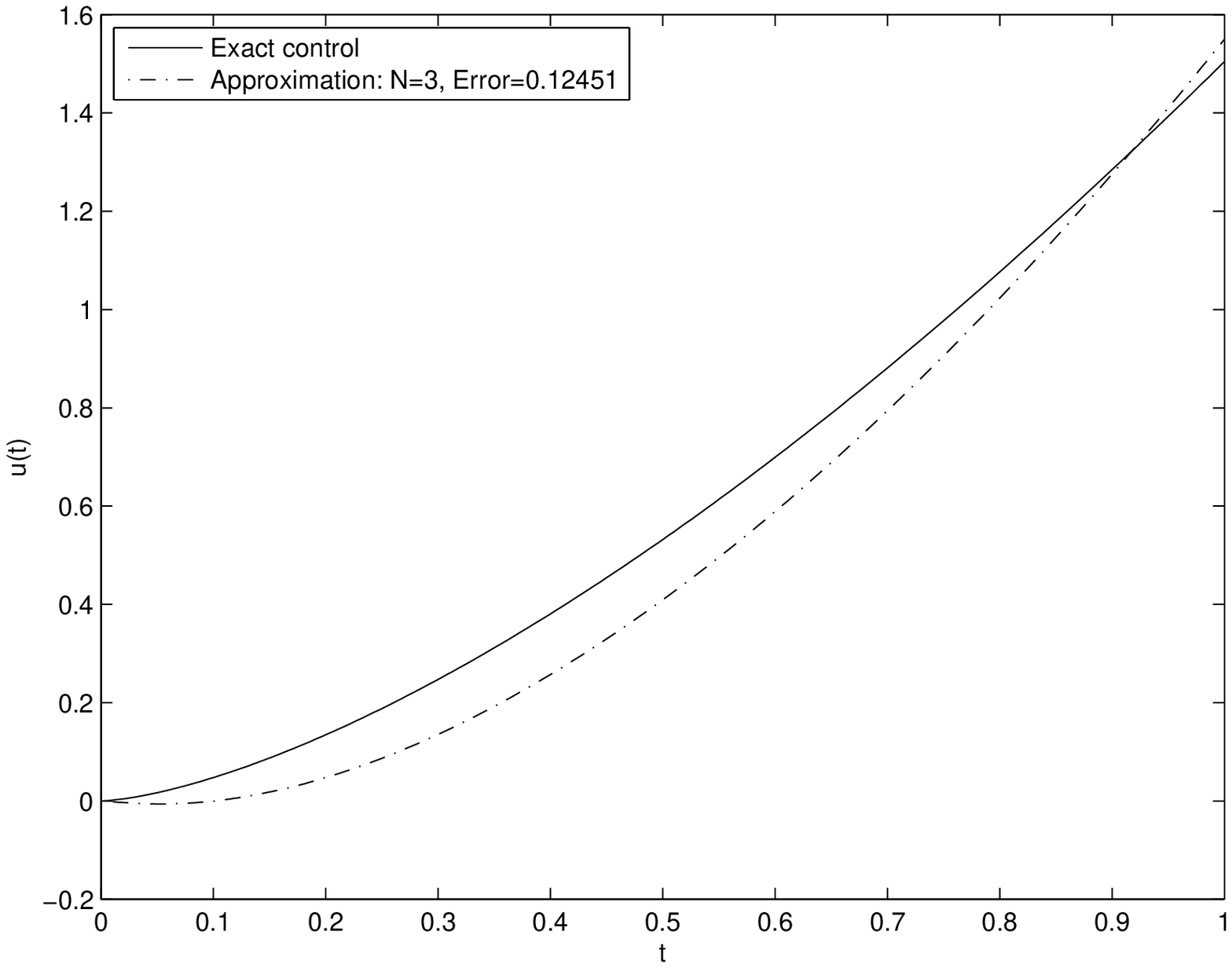}}
\end{center}
\caption{\label{plot5} Exact solution (solid lines) for the fractional 
optimal control problem \eqref{eq:N1}--\eqref{eq:N3} with $\alpha = 1/2$
versus numerical solutions (dashed lines) obtained approximating the
optimality conditions given by Theorem~\ref{teo1}.}
\end{figure}

We will see that applying the numerical method given in the previous section,
we are able to solve fractional optimal control problems applying
known techniques from the classical optimal control theory.
First, consider the following optimal control problem:
\begin{equation}
\label{eq:N1}
J(x,u)=\int_0^1 \left(t u(t)-(\a+2)x(t)\right)^2\,dt \longrightarrow \min
\end{equation}
subject to the control system
\begin{equation}
\label{eq:N2}
\dot{x}(t)+{^C_0D^\a_t} x(t)=u(t)+t^2
\end{equation}
and the boundary conditions
\begin{equation}
\label{eq:N3}
x(0)=0, \quad x(1)=\frac{2}{\Gamma(3+\a)}.
\end{equation}
The solution is given by
$$
\left(\overline x(t),\overline u(t)\right)
=\left(\frac{2t^{\a+2}}{\Gamma(\a+3)},\frac{2t^{\a+1}}{\Gamma(\a+2)}\right).
$$
Using the necessary conditions given in Theorem~\ref{teo1}, we arrive at
$$
\begin{cases}
\dot{x}(t)+{^C_0D^\a_t} x(t)=-\frac{\l}{2t^2}+\frac{\a+2}{t}x(t)+t^2,\\
-\dot{\l}(t)+{_tD^\a_1}\l(t)=\frac{\a+2}{t}\l(t),
\end{cases}
\quad
\begin{cases}
x(0)=0,\\
x(1)=\frac{2}{\Gamma(3+\a)},
\end{cases}
$$
which is a fractional boundary value problem.
We approximate this problem by approximation 
in \eqref{2} up to order $N$:
$$
\begin{cases}
\dot{x}(t) = \left[\left(\frac{\a+2}{t}-At^{-\a}\right)x(t)
+\sum_{p=2}^N C_pt^{1-p-\a}V_p(t)-\frac{\l(t)}{2t^2}
+t^2\right] \times\frac{1}{1+Bt^{1-\a}}\\
\dot{V}_p(t)= (1-p)t^{p-2}x(t), \quad p=2,\ldots,N\\
\dot{\l}(t)=\left[\left(A(1-t)^{-\a}-\frac{\a+2}{t}\right)\l(t)
-\sum_{p=2}^N C_p(1-t)^{1-p-\a}W_p(t)\right]\times \frac{1}{1+B(1-t)^{1-\a}}\\
\dot{W}_p(t)= -(1-p)(1-t)^{p-2}\l(t), \quad p=2,\ldots,N
\end{cases}
$$
subject to the boundary conditions
$$
\begin{cases}
x(0)=0,\quad x(1)=\frac{2}{\Gamma(3+\a)},\\
V_p(0)=0, \quad p=2,\ldots,N,\\
W_p(1)=0, \quad p=2,\ldots,N.
\end{cases}
$$
The solutions are depicted in Figure~\ref{plot5},
for $N=2$, $N=3$, and $\a=1/2$, with the error being given by
$E = \max_{i}(|x(t_i)-\overline x(t_i)|)$.
\begin{figure}[ht]
\begin{center}
\subfigure[$x(t), N=2$]{\includegraphics[scale=0.37]{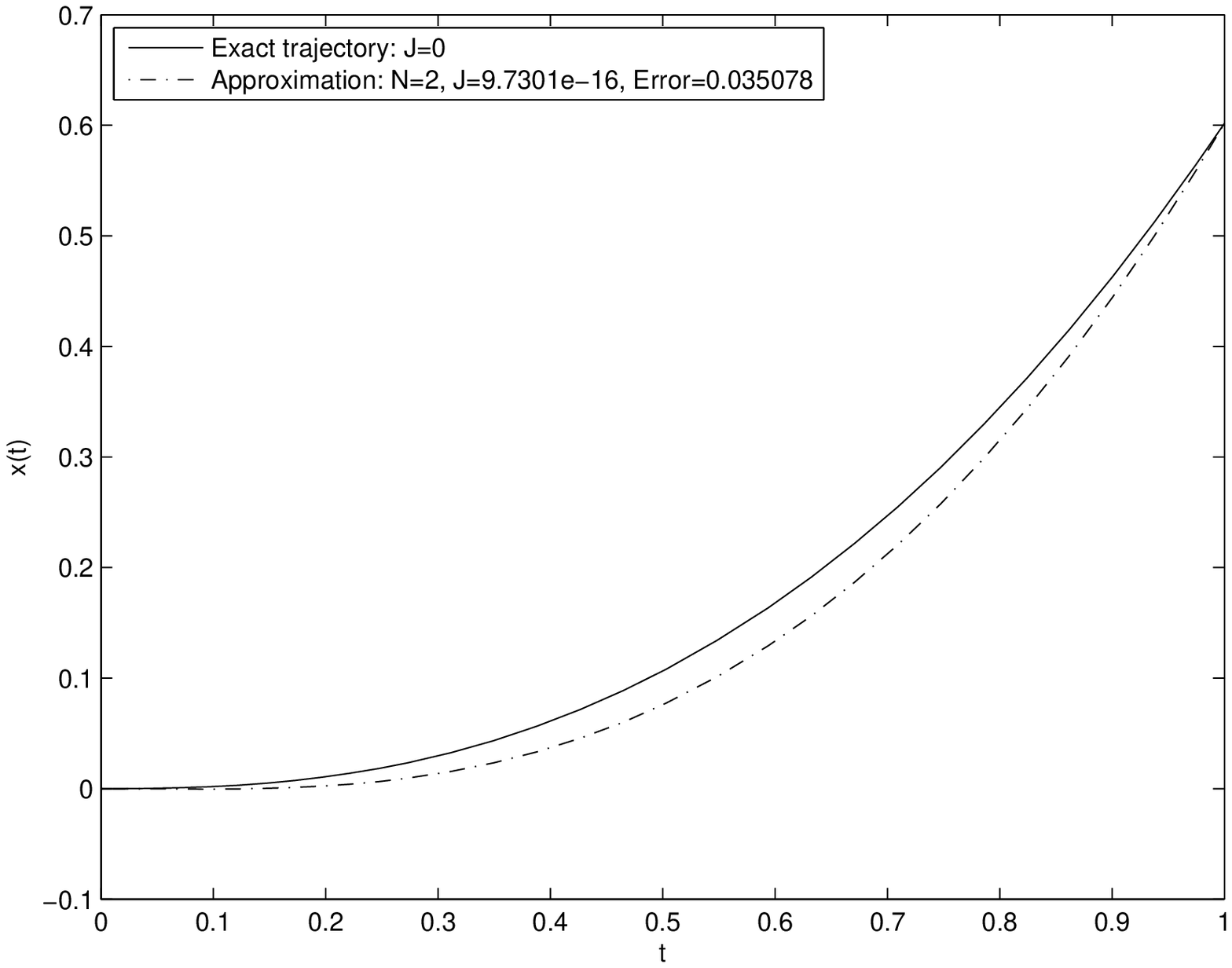}}
\subfigure[$u(t), N=2$]{\includegraphics[scale=0.37]{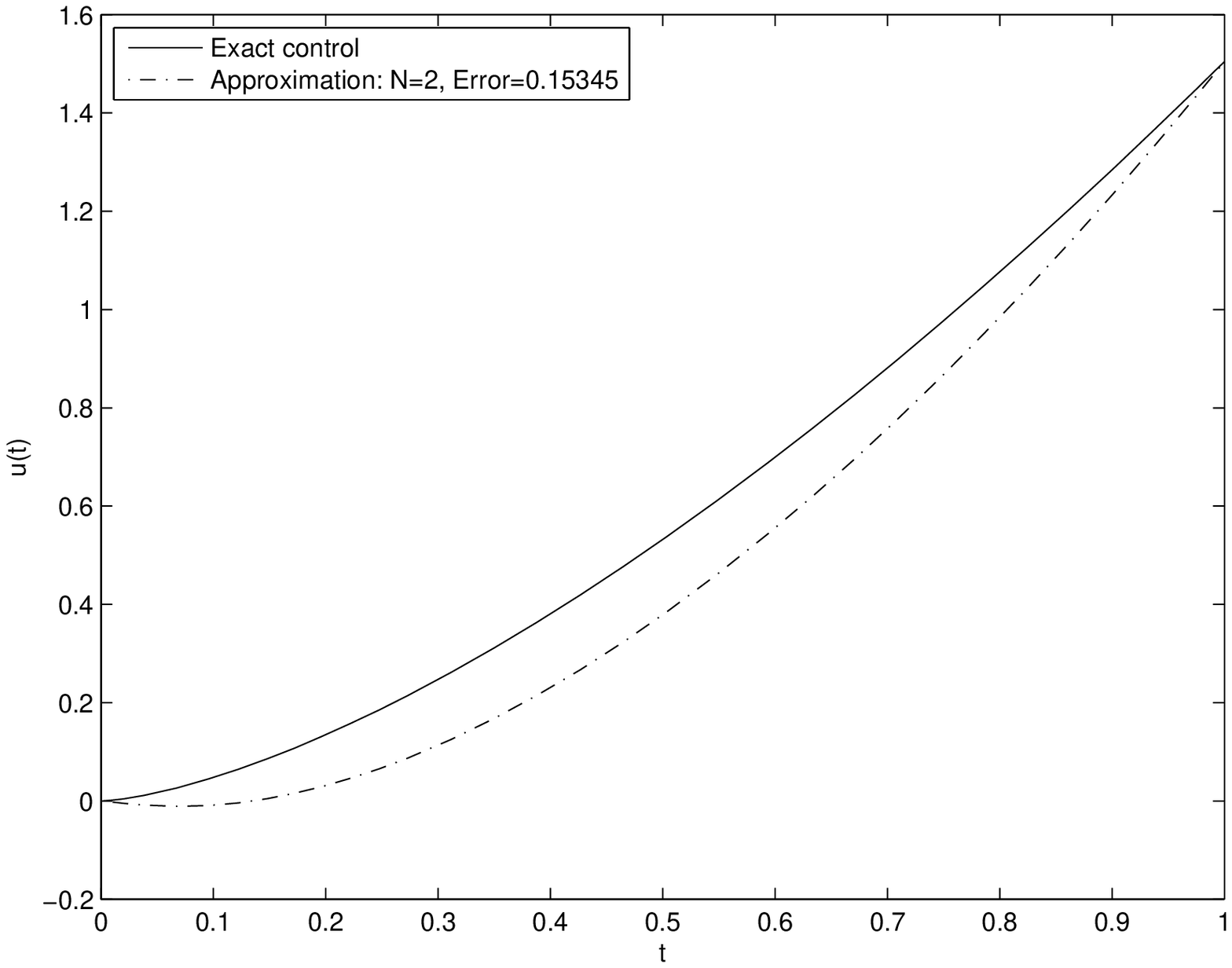}}
\end{center}
\end{figure}

\begin{figure}[ht]
\begin{center}
\subfigure[$x(t), N=2$]{\includegraphics[scale=0.5]{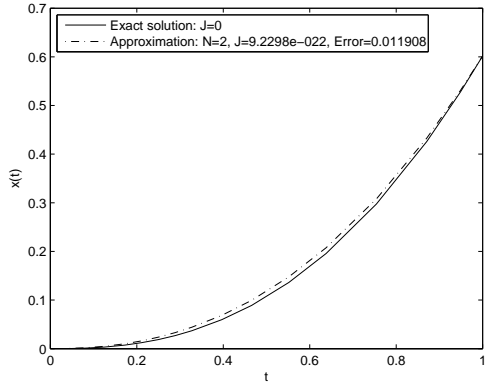}}
\subfigure[$u(t), N=2$]{\includegraphics[scale=0.5]{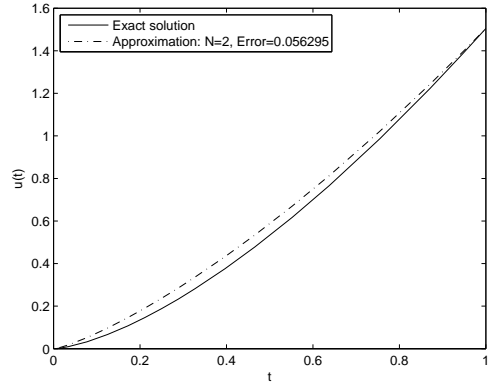}}
\subfigure[$x(t), N=3$]{\includegraphics[scale=0.5]{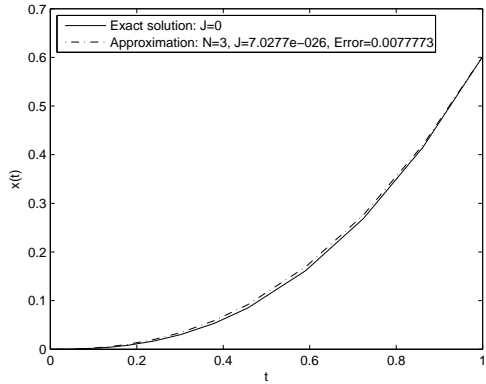}}
\subfigure[$u(t), N=3$]{\includegraphics[scale=0.5]{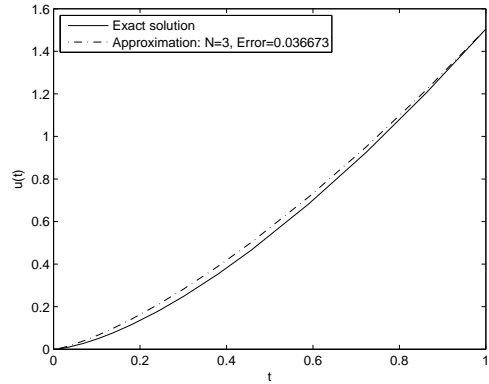}}
\end{center}
\caption{\label{plot6}Exact solution (solid lines) 
for the fractional optimal control problem \eqref{eq:N1}--\eqref{eq:N3} 
with $\alpha = 1/2$ versus numerical solutions (dashed lines) obtained
approximating the original problem by a classical one.}
\end{figure}

Another approach is to approximate the original problem
by using approximation from \eqref{2} directly, getting
$$
\tilde{J}(x,u)=\int_0^1(tu-(\a+2)x)^2\,dt \longrightarrow \min
$$
subject to the control system
$$
\begin{cases}
\dot{x}(t)[1+B(\a,N)t^{1-\a}]+A(\a,N)t^{-\a}x(t)
-\sum_{p=2}^N C(\a,p)t^{1-p-\a}V_p(t)=u(t)+t^2\\
\dot{V}_p(t)=(1-p)t^{p-2}x(t)
\end{cases}
$$
and boundary conditions
$$
x(0)=0, \quad x(1)=\frac{2}{\Gamma(3+\a)},
\quad V_p(0)=0, \quad p=2,3,\ldots,N.
$$
The (classical) necessary optimality conditions become
$$
\begin{cases}
\dot{x}(t)=2\phi_0(t)\l_1(t)+\phi_1(t)x(t)
+\sum_{p=2}^N \phi_p(t)V_p(t)+\phi_{N+1}(t)\\
\dot{V_p}=(1-p)t^{p-2}x(t), \quad p=2,\ldots,N\\
\dot{\lambda}_1=-\phi_1(t)\l_1(t)+\sum_{p=2}^N(p-1)t^{p-2}\l_p\\
\dot{\lambda}_p=-\phi_p(t)\l_1(t), \quad p=2,\ldots,N
\end{cases}
$$
subject to the boundary conditions
$$
\begin{cases}
x(0)=0\\
V_p(0)=0, \quad p=2,\ldots,N
\end{cases}
\qquad
\begin{cases}
x(1)=\frac{2}{\Gamma(3+\a)}\\
\l_p(1)=0, \quad p=2,\ldots,N.
\end{cases}
$$
The solutions are depicted in Figure~\ref{plot6} 
for $N=2$, $N=3$, and $\a=1/2$.

For our next example, we consider the final time $T$ free and thus a variable in the problem.
We wish to find an optimal triplet $(x(\cdot),u(\cdot),T)$ that minimizes
\begin{equation}
\label{eq:N4}
J(x,u,T)=\int_0^T(tu-(\a+2)x)^2\,dt
\end{equation}
subject to the control system
\begin{equation}
\label{eq:N5}
\dot{x}(t)+{^C_0D^\a_t} x(t)=u(t)+t^2
\end{equation}
and boundary conditions
\begin{equation}
\label{eq:N6}
x(0)=0, \quad x(T)=1.
\end{equation}
In this case, an exact solution to this problem is not known.

The fractional necessary optimality conditions,
after approximating the fractional terms, result in
$$
\begin{cases}
\dot{x}(t)=\left[\left(\frac{\a+2}{t}-At^{-\a}\right)x(t)
+\sum_{p=2}^N C_pt^{1-p-\a}V_p(t)
-\frac{\l(t)}{2t^2}+t^2\right]\times \frac{1}{1+Bt^{1-\a}}\\
\dot{V}_p(t)=(1-p)t^{p-2}x(t), \quad p=2,\ldots,N\\
\dot{\l}(t)=\left[\left(A(1-t)^{-\a}-\frac{\a+2}{t}\right)\l(t)
-\sum_{p=2}^N C_p(1-t)^{1-p-\a}W_p(t)\right] \times\frac{1}{1+B(1-t)^{1-\a}}\\
\dot{W}_p(t)=-(1-p)(1-t)^{p-2}\l(t), \quad p=2,\ldots,N
\end{cases}
$$
subject to the boundary conditions
$$
\begin{cases}
x(0)=0,\quad x(T)=1,\\
V_p(0)=0, \quad p=2,\ldots,N,\\
W_p(T)=0, \quad p=2,\ldots,N.
\end{cases}
$$

\begin{figure}[ht]
\begin{center}
\subfigure[$x(t), N=2$]{\includegraphics[scale=0.37]{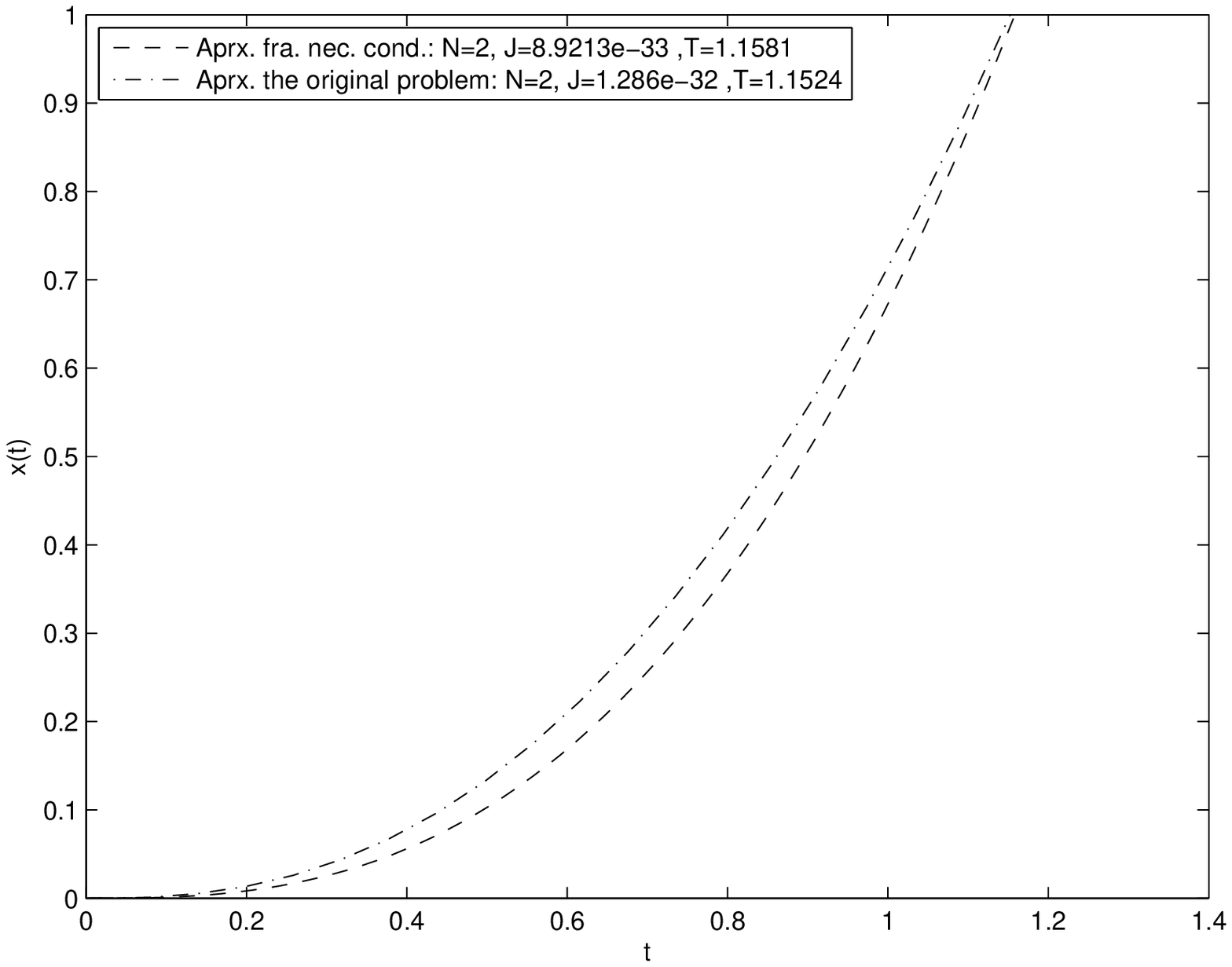}}
\subfigure[$u(t), N=2$]{\includegraphics[scale=0.37]{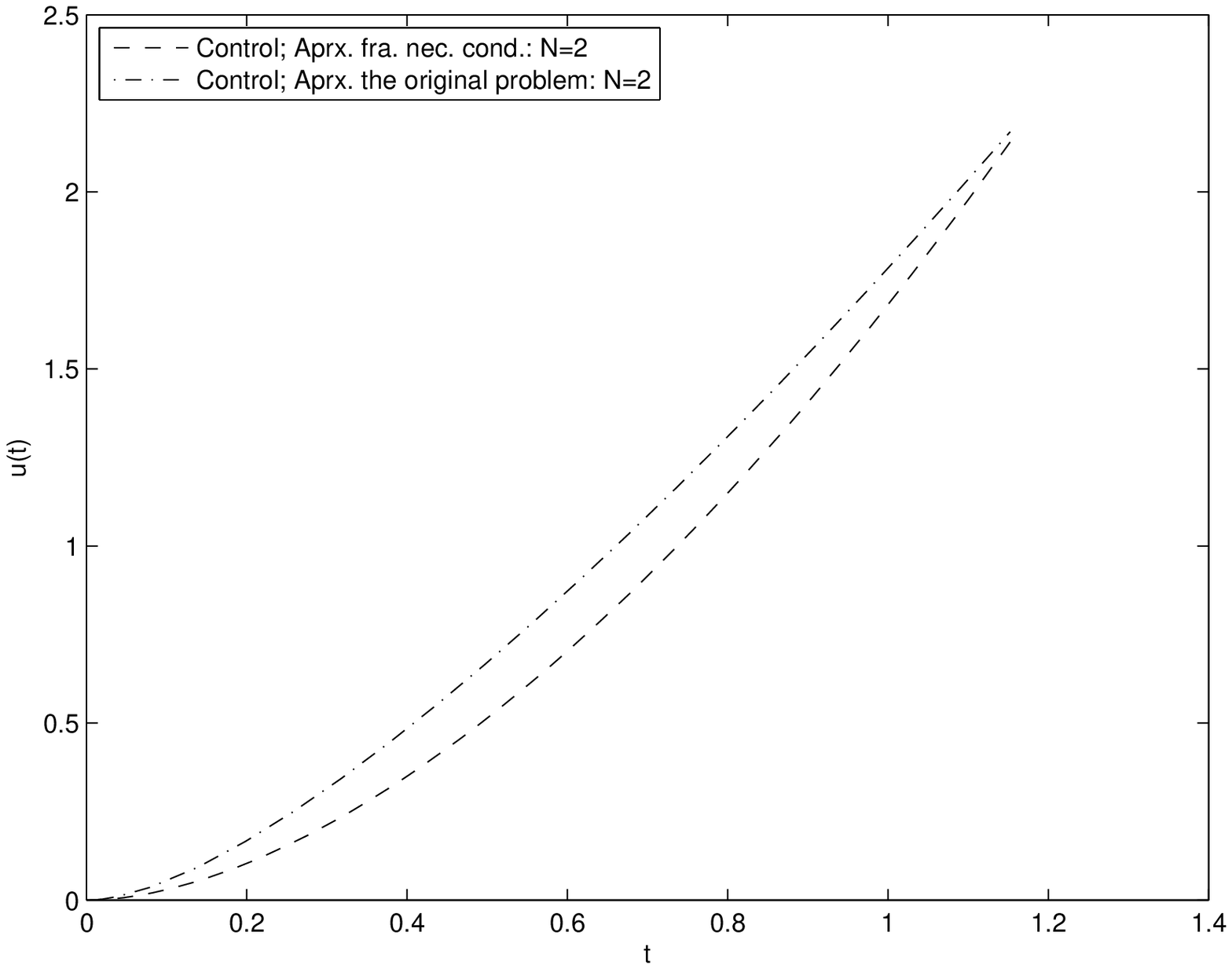}}
\end{center}
\caption{\label{plot7}Numerical solutions to the free final time problem 
\eqref{eq:N4}--\eqref{eq:N6}, using fractional necessary optimality conditions 
(dashed lines) and approximation of the problem to an integer order 
optimal control problem (dash-dotted lines).}
\end{figure}

Another way is transforming the problem into an integer order optimal
control problem with free final time. The necessary optimality conditions are
$$
\begin{cases}
\dot{x}(t)=2\phi_0(t)\l_1(t)+\phi_1(t)x(t)+\phi_2(t)V_2(t)+\phi_{3}(t)\\
\dot{V_2}=-x(t)\\
\dot{\lambda}_1=-\phi_1(t)\l_1(t)+x(t)\\
\dot{\lambda}_2=-\phi_2(t)\l_1(t).
\end{cases}
$$
The results obtained are shown in Figure~\ref{plot7}.


\section*{Acknowledgments}

This work was supported by {\it FEDER} funds through {\it COMPETE}
--- Operational Programme Factors of Competitiveness
(``Programa Operacional Factores de Competitividade'')
and by Portuguese funds through the {\it Center for Research
and Development in Mathematics and Applications} (University of Aveiro)
and the Portuguese Foundation for Science and Technology
(``FCT--Funda\c{c}\~{a}o para a Ci\^{e}ncia e a Tecnologia''),
within project PEst-C/MAT/UI4106/2011 with COMPETE
number FCOMP-01-0124-FEDER-022690. Pooseh was also
supported by the FCT Ph.D. fellowship SFRH/BD/33761/2009.




\begin{thebibliography}{99}

\bibitem{kilbas}
A. A. Kilbas, H. M. Srivastava\ and\ J. J. Trujillo,
{\it Theory and applications of fractional differential equations},
North-Holland Mathematics Studies, 204, Elsevier, Amsterdam, 2006.

\bibitem{Podlubny}
I. Podlubny,
{\it Fractional differential equations},
Mathematics in Science and Engineering,
198, Academic Press, San Diego, CA, 1999.

\bibitem{b:agnieszka:delfim}
A. B. Malinowska\ and\ D. F. M. Torres,
{\it Introduction to the fractional calculus of variations},
Imp. Coll. Press, London, 2012.

\bibitem{Atanackovic}
T. M. Atanackovic\ and\ B. Stankovic,
On a numerical scheme for solving differential equations of fractional order,
Mech. Res. Comm. {\bf 35} (2008), no.~7, 429--438.

\bibitem{Almeida1}
S. Pooseh, R. Almeida\ and\ D. F. M. Torres,
Expansion formulas in terms of integer-order derivatives
for the Hadamard fractional integral and derivative,
Numer. Funct. Anal. Optim. {\bf 33} (2012), no.~3, 301--319.
{\tt arXiv:1112.0693}

\bibitem{Almeida3}
S. Pooseh, R. Almeida\ and\ D. F. M. Torres,
Approximation of fractional integrals by means of derivatives,
Comput. Math. Appl. {\bf 64} (2012), no.~10, 3090--3100.
{\tt arXiv:1201.5224}

\bibitem{Almeida2}
S. Pooseh, R. Almeida\ and\ D. F. M. Torres,
Numerical approximations of fractional derivatives with applications,
Asian J. Control {\bf 15} (2013), no.~3, 698--712.
{\tt arXiv:1208.2588}

\bibitem{Almeida4}
S. Pooseh, R. Almeida\ and\ D. F. M. Torres,
Fractional order optimal control problems with free terminal time,
J. Ind. Manag. Optim., in press.
{\tt arXiv:1302.1717}

\bibitem{Almeida0}
R. Almeida, S. Pooseh\ and\ D. F. M. Torres,
Fractional variational problems depending on indefinite integrals,
Nonlinear Anal. {\bf 75} (2012), no.~3, 1009--1025.
{\tt arXiv:1102.3360}

\bibitem{Almeida5}
S. Pooseh, R. Almeida\ and\ D. F. M. Torres,
Discrete direct methods in the fractional calculus of variations,
Comput. Math. Appl., in press. DOI:10.1016/j.camwa.2013.01.045
{\tt arXiv:1205.4843}

\end{thebibliography}
\end{document}